\documentclass[12pt,reqno]{amsart}

\usepackage{amssymb}

\newcommand{\ds}{\displaystyle}     

\DeclareMathAlphabet{\E}{U}{eus}{m}{n}     

\newcommand{\V}{{\mathcal V}}

\newcommand{\PP}{{\mathbb P}}

\newcommand{\kk}{{\Bbbk}}

\renewcommand{\char}{\mbox{\rm char}}

\newtheorem{thm}{Theorem}[section]
\newtheorem{lemma}[thm]{Lemma}
\newtheorem{cor}[thm]{Corollary}
\newtheorem{prop}[thm]{Proposition}

\theoremstyle{definition}
\newtheorem{defn}[thm]{Definition}

\newtheorem{ex}[thm]{Example}
\newtheorem{rmk}[thm]{Remark}

\newtheorem{ques}[thm]{Question}

\newcommand{\dual}{\raisebox{0.4ex}{$^*$}\!\!}
\newcommand{\sstar}{\raisebox{0.4ex}{$\star$}}

\newtheorem{Pf}{Proof$\!\!$}         
         
\newenvironment{pf}{\begin{Pf}}{\qed\end{Pf}}

\DeclareMathSymbol{\twoheadrightarrow}  {\mathrel}{AMSa}{"10}

\newcounter{letter}
\renewcommand{\theletter}{\rom{(}\alph{letter}\rom{)}}

\newcounter{rnum}
\renewcommand{\thernum}{\rom{(}\roman{rnum}\rom{)}}

\begin{document}
\baselineskip21pt

\title[Quantum $\PP^2$s using Graded Skew Clifford Algebras]%
{Classifying Quadratic Quantum $\PP^2$s by using\\[3mm] Graded Skew 
Clifford Algebras}


\maketitle

\vspace*{0.1in}

\baselineskip18pt

\centerline{\sc Manizheh Nafari\footnote{M.~Nafari was
supported in part by NSF grant DMS-0457022.} 
\quad {\sf manizheh@uta.edu}}
\centerline{\sc Michaela Vancliff\footnote{M.~Vancliff was
supported in part by NSF grants DMS-0457022 \& DMS-0900239.}
\quad {\sf vancliff@uta.edu \quad uta.edu/math/vancliff}}
\centerline{\sc Jun Zhang \quad {\sf zhangjun19@gmail.com}}
\quad\\[0mm]
\centerline{Department of Mathematics, P.O.~Box 19408,}
\centerline{University of Texas at Arlington,}
\centerline{Arlington, TX \ 76019-0408}

\setcounter{page}{1}
\thispagestyle{empty}

\bigskip
\bigskip

\begin{abstract}
\baselineskip15pt
We prove that quadratic regular algebras of global dimension three on 
degree-one generators are related to graded skew Clifford algebras. In
particular, we prove that almost all such algebras may be constructed as a 
twist of either a regular graded skew Clifford algebra or of an Ore extension
of a regular graded skew Clifford algebra of global dimension two.
In so doing, we classify all quadratic regular algebras of 
global dimension three that have point scheme either a nodal cubic curve or 
a cuspidal cubic curve in $\PP^2$.
\end{abstract}

\baselineskip24pt


\bigskip
\bigskip

\section*{Introduction}

In \cite{AS}, M.~Artin and W.~Schelter introduced a notion of regularity 
for a non-commutative graded algebra on degree-one generators. To such an 
algebra one may associate some geometry via certain graded modules over the 
algebra, as discussed by M.~Artin, J.~Tate and M.~Van den Bergh in 
\cite{ATV1}. In this spirit, one may describe such an algebra using a certain
scheme (called the point scheme in \cite{VV})
and, in \cite{A}, M.~Artin introduced the language 
of ``quantum $\PP^2$'' for a regular algebra of global dimension three 
on degree-one generators. Generic quantum $\PP^2$s were classified in 
\cite{AS,ATV1,ATV2}. 

In \cite{CV}, a new relatively simple method was given for constructing 
some quadratic regular algebras of finite global dimension, and
quadratic regular algebras 
produced with that method are called regular graded skew Clifford algebras
(see Definition~\ref{GSCA}).  At this time, it is unclear 
how useful this method will be in helping resolve the open problem of 
classifying all quadratic regular algebras of global dimension four. 
A first step in this direction is to determine how useful graded skew 
Clifford algebras are in classifying quadratic quantum~$\PP^2$s.   In
this article, we prove that most quadratic quantum~$\PP^2$s
may be constructed as a twist of 
either a regular graded skew Clifford algebra or of an Ore extension of
a regular graded skew Clifford algebra of global dimension two.
For the precise results,
see Theorem~\ref{main1'}, Corollary~\ref{cor2}, Corollary~\ref{cor3}, 
Lemma~\ref{H}, Lemma~\ref{B} and Lemma~\ref{A}.  
The only algebras we did not relate in this way to graded skew Clifford 
algebras are some that have point scheme an elliptic curve; specifically 
those of type~E in \cite{AS} (up to isomorphism and anti-isomorphism, there 
is at most one such algebra) and an open subset of those of type~A in 
\cite{AS} (although a weak relationship is discussed for these algebras
in Remark~\ref{weak}).
It is still open whether or not similar results hold for quadratic regular 
algebras of global dimension four.

Henceforth, let $A$ denote a quadratic quantum $\PP^2$.
The work in this article is partitioned across sections according to the 
point 
scheme of $A$.  Those $A$ whose point scheme contains a line are discussed 
in section~\ref{Sec1}; is a nodal cubic curve in section~\ref{Sec2}; cuspidal 
cubic curve in section~\ref{Sec3}; and, finally, elliptic curve in 
section~\ref{Sec4}.   
Those $A$ whose point scheme is either a nodal cubic curve or a cuspidal 
cubic curve in $\PP^2$ are not specifically discussed in 
\cite{ATV1}, as they are not generic. In Sections 2 and 3 of this article, we 
classify all such algebras (see Theorem~\ref{nodal} and Theorem~\ref{cusp}). 
Indeed, up to isomorphism, there is at most a one-parameter family of 
quadratic quantum $\PP^2$s with point scheme a nodal cubic curve, whereas 
there is only one such algebra with point scheme a cuspidal cubic curve.


\bigskip
\bigskip

\section{Quantum $\PP^2$s with Reducible or Non-reduced Point Scheme}%
\label{Sec1}

Throughout the article, $\kk$ denotes an algebraically closed field and 
$\kk^\times = \kk\setminus \{ 0 \}$. 
It is well known that a quadratic regular algebra of global dimension $n+1$
with point scheme given by $\PP^n$ is a twist, by an automorphism, of the 
polynomial ring on $n+1$ variables (c.f., \cite[Page 378]{ATV2} and 
Remark~\ref{rmk1}(iii)). 
We prove a generalization of this result in this section for certain 
quadratic quantum $\PP^2$s. Precisely, we prove (Theorem~\ref{main1}) that 
a quadratic quantum $\PP^2$ that has a reducible 
or non-reduced point scheme is either a twist by an automorphism of a
skew polynomial ring or is a twist, by a twisting system, of an Ore extension
of a regular algebra of global dimension two. This result is restated in 
terms of graded skew Clifford algebras in Theorem~\ref{main1'} where 
char$(\kk) \neq 2$.

As in the Introduction, let $A$ denote a quadratic quantum $\PP^2$ with 
point scheme given by a cubic divisor $C \subset \PP^2$.  (Here, $\PP^2$
may be identified with $\PP(A\dual_1)$ where $A_1$ denotes the span of the 
homogeneous degree-one elements of $A$, and $A\dual_1$ denotes the 
vector-space dual of $A_1$.) Throughout this section, we assume that $C$ is 
reducible or non-reduced, so that $C$ is either the union of a nondegenerate 
conic and a line, or the union of three distinct lines, or the union of a 
double line and a line, or is a triple line. The automorphism encoded 
by the point scheme will be denoted by $\sigma \in$ Aut$(C)$. There are
two cases to consider: either $C$ contains a line that is invariant
under $\sigma$ or it does not.  Both cases use the notion of twisting a
graded algebra by an automorphism, which is defined in \cite[\S8]{ATV2};
in the case of a quadratic algebra, it is defined as follows.

\begin{defn}\label{twistdef}\cite[\S8]{ATV2} \
Let $D$ denote a quadratic algebra, let $D_1$ denote the span of the
homogeneous degree-one elements of $D$ and let $\phi$ be a graded degree-zero 
automorphism of~$D$.  The twist $D^\phi$ of $D$ by $\phi$ is a quadratic
algebra that has the same underlying vector space as~$D$, but has a new 
multiplication $*$ defined as follows: if $a, b \in D_1 = (D^\phi)_1$, then 
$a * b = a \phi(b)$, where the right-hand side is computed using the original 
multiplication in $D$.
\end{defn}
In general, twisting by an automorphism is a reflexive and symmetric
operation, but not a transitive operation; in fact, 
twisting an algebra~$D$ twice yields a twist of $D$ by a 
{\em twisting system}, and that notion is defined in \cite{Z}.

\subsection{Case 1: $C$ contains a line invariant under $\sigma$}\hfill

Suppose $L\subset C$ is a line that is invariant under $\sigma$, and let
$x \in A_1$ be such that $L$ is the zero locus, $\V(x)$, of $x$.   By 
\cite[Theorem 8.16(i)(ii) and Corollary 8.6]{ATV2}, $x$ is normal in $A$, and
one may twist $A$, by an automorphism, to obtain a quadratic regular algebra 
$B$ in which the image of $x$ is central.

\begin{prop}\label{prop1}
In the above notation, the twist $B$ of $A$ is an Ore extension of the
polynomial ring on two variables.
\end{prop}
\begin{pf}
Let $x'$ denote the image of $x$ in $B$.  Since $B/\langle x' \rangle$ is a 
regular algebra of global dimension two, it is isomorphic to either 
$\kk\langle Y, Z\rangle/\langle YZ - ZY - Y^2\rangle$ or 
$\kk\langle Y, Z\rangle/\langle ZY - qYZ\rangle$ where $q \in\kk^\times$
(\cite[Page 172]{AS}). It follows from \cite[Theorem 8.16(iii)]{ATV2} that 
$B$ is a $\kk$-algebra on generators $x', y, z$ with defining 
relations\\[-7mm]
\[
x'y = yx', \qquad x'z = zx', \qquad h = 0,
\]
where either 
\[
\begin{array}{rl}
\text{(i)}& h = zy - yz + y^2 + x'(\alpha x' + \beta y + \gamma z)
\text{\ \ or}\\[3mm]
\text{(ii)}& h = zy - qyz + x'(\alpha x' + \beta y + \gamma z),\\[2mm]
\end{array}
\]
for some $\alpha$, 
$\beta$, $\gamma \in \kk$. Let $B' = \kk[x' , y]$, and define $\phi \in$
Aut$(B')$ and a left $\phi$-derivation $\delta:B'\to B'$ as follows:
\[
\phi (x') = x',\qquad \phi (y) = s y - \gamma x', \qquad \delta (x') = 0,%
\qquad \delta (y) = -t y^2 -\alpha (x')^2 - \beta x'y,
\] 
where $(s,\ t) = (1,\ 1)$ if $h$ is given by (i) and 
$(s,\ t) = (q ,\ 0)$ if $h$ is given by (ii). 
In both cases, $B = B'[z; \phi, \delta]$.
\end{pf}

\begin{cor}\label{cor1}
If the point scheme of $A$ contains a line that is invariant under $\sigma$,
then $A$ is a twist, by an automorphism, of an Ore extension of the
polynomial ring on two variables. 
\end{cor}
\begin{pf}
Combining Proposition~\ref{prop1} with the preceding discussion proves 
the result.
\end{pf}

\subsection{Case 2: no line in $C$ is invariant under $\sigma$}
\hfill

Suppose that no line in $C$ is invariant under $\sigma$. It follows that
$C$ is the union of three distinct lines that are cyclically permuted by
$\sigma$.  Such lines have the property that either no point lies on all
three lines or the three lines meet at exactly one point.

\begin{rmk} \label{rmk1}\hfill\\
\indent (i)
Let $D$ denote a quadratic quantum $\PP^2$.  Let $V = D_1$ and let 
$W \subset V \otimes_\kk V$ denote the span of the defining relations of 
$D$, and let $\V(W)$ denote the zero locus, in $\PP(V\dual\ ) \times 
\PP(V\dual\ )$, 
of the elements of $W$, where $\PP(V\dual\ )$ is identified with $\PP^2$.
The Koszul dual $D^{\sstar}$ of $D$ is the quotient 
of the free algebra on $V\dual\ $ by the ideal generated by $W^\perp
\subset V\dual\ 
\otimes_\kk V\dual\ $.  If $(p, q) \in \V(W)$, then $\kk p \subset
V\dual\ $,
$\kk q \subset V\dual\ $ and $p \otimes q \in W^\perp$, and conversely.  This 
provides a method of passing between $\V(W)$ and the relations of 
$D^{\sstar}$.

(ii) With notation as in (i), by \cite{ATV1}, $\V(W)$ is the graph of 
an automorphism $\tau$ of a subscheme $\mathcal P$ of $\PP(V\dual\ ) = \PP^2$. 

(iii) With notation as in (i) and (ii), if $\tau$ may be extended to an 
automorphism of $\PP(V\dual\ ) = \PP^2$, then, since $D$ is regular, $D$ is a 
twist, by an automorphism, of the polynomial ring and $\mathcal 
P = \PP(V\dual\ )$. This is because the homogeneous degree-two forms that vanish 
on the graph of $\tau$ have the form $\tau(u) v - \tau(v) u$ for all $u,\ v 
\in D_1$. 
\end{rmk}

\begin{lemma}\label{lemma1}
Suppose $C$ is the union of three distinct lines that are cyclically 
permuted by $\sigma$.  If no point lies on all three lines, 
then $A$ is a twist, by an automorphism, of a skew polynomial ring.
\end{lemma}
\begin{pf}
By hypothesis, there exist linearly independent elements  
$x, y, z \in A_1$ such that $C = \V(xyz)$ and $\sigma : 
\V(x) \to \V(y) \to \V(z) \to \V(x)$. Since $\sigma(1, 0, 0) \in \V(x, z)$
and $\sigma(0, 1, 0) \in \V(x, y)$ and $\sigma(0, 0, 1) \in \V(y, z)$, it
follows that 
\[
\sigma(0, \beta, \gamma) = (d \gamma , 0, \beta),\quad
\sigma(\alpha, 0, \gamma) = (\gamma , e \alpha, 0),\quad
\sigma(\alpha, \beta, 0) = (0 , \alpha, f \beta),
\]
for some $d, e, f \in \kk^\times$, for all $(\beta, \ \gamma)$, 
$(\alpha, \ \gamma)$, $(\alpha, \ \beta) \in \PP^1$. 
By Remark~\ref{rmk1}(ii), this implies that $A$ is a $\kk$-algebra
on generators $x, y, z$ with defining relations:
\[
yx = d z^2, \qquad zy = e x^2, \qquad xz = f y^2.
\]
Define $\tau \in$ Aut$(A)$ by $\tau(x) = \lambda_1 y$, $\tau (y)=\lambda_2 z$
and $\tau(z) = \lambda_3 x$, where $\lambda_1$, $\lambda_2$ $\lambda_3\in
\kk^\times$ satisfy $d \lambda_1\lambda_3^2 = e \lambda_1^2\lambda_2= 
f \lambda_2^2 \lambda_3$.  Twisting by $\tau$ yields a $\kk$-algebra
on $x, y, z$ with defining relations:
\[
yz + c zy =0, \qquad zx + c xz=0, \qquad xy + c yx=0,
\]
for some $c \in \kk^\times$, and such an algebra is a skew polynomial ring.
\end{pf}

\begin{prop}\label{prop2}
Suppose $C$ is the union of three distinct lines that are cyclically 
permuted by $\sigma$.  If all three lines intersect at exactly one point,
then $A$ is either\/ {\rm (a)} an Ore extension of a regular algebra of global
dimension two, or\/ {\rm (b)} a twist, by a twisting system, of an Ore 
extension of the polynomial ring on two variables; if\/ $\char(\kk) \neq 3$, 
then $A$ is described by\/ {\rm (b)}.
\end{prop}
\begin{pf}
Suppose all three lines of $C$ intersect at only one point $p=(0, 0, 1)$,
and write $C = \V(xy(x-y))$, where $x, y \in A_1$. We may assume 
$\sigma : \V(x) \to \V(y) \to \V(x-y) \to \V(x)$.  Since $\sigma(p) = p$, 
we have 
\[
\sigma(0, \beta, \gamma) = (\beta , 0, a\beta + b \gamma),\quad
\sigma(\alpha, 0, \gamma) = (\alpha , \alpha, c \alpha + d \gamma),\quad
\sigma(\beta, \beta, \gamma) = (0 , \beta, e \beta + f \gamma),
\]
for some $a, c, e \in \kk$, $b, d, f \in \kk^\times$, 
for all $(\beta, \ \gamma)$, $(\alpha, \ \gamma) \in \PP^1$. 
In $A\dual_1$, let $\{ X, Y, Z\}$ denote the dual basis to $\{x, y, z\}$.
By Remark~\ref{rmk1}(i), the following relations hold in the Koszul dual of 
$A$: 
\begin{gather*}
YX + a YZ = 0, \quad Y^2 + XY + eXZ + e YZ= 0, \quad X^2 + XY + cXZ = 0,\\
Z^2 = 0, \quad 
ZY + ZX + d XZ= 0, \quad ZX + b YZ = 0, \quad ZY + fYZ + f XZ= 0.
\end{gather*}
Since $A$ is a quantum $\PP^2$, these relations span at most a 6-dimensional 
space. As the first four relations span a 4-dimensional space, and the last 
three relations are linearly independent of the first four relations, it 
follows 
that the span of the last three relations has dimension at most two. This 
implies that $d = f = -b$, so that we may write
\[
\sigma(0, \beta, \gamma) = (\beta , 0, a\beta + b \gamma),\quad
\sigma(\alpha, 0, \gamma) = (\alpha , \alpha, c \alpha - b \gamma),\quad
\sigma(\beta, \beta, \gamma) = (0 , \beta, e \beta - b \gamma).
\]
By Remark~\ref{rmk1}(iii), since $C \neq \PP^2$, $\sigma$ cannot be extended 
to $\PP^2$, from which we obtain $c \neq e + a$.
Suppose $b^2 + b + 1 \neq 0$, and define
\[\tau = \begin{bmatrix} 0 & -1 & 0\\ 1 & -1 & 0 \\ k & g & 1 \end{bmatrix}
\in \text{Aut}(C),\]
for some $g, k \in \kk$.  Thus,
$\tau : \V(x) \to \V(x-y) \to \V(y) \to \V(x)$.  
In order for $\sigma$ and $\tau$ to commute on $C$, we require $g$ and $k$
to satisfy
\[
bg - k = a - e \quad \text{and} \quad g + (b+1)k = -a -c.
\]
Since $b^2 + b + 1 \neq 0$, these equations have a unique solution, so we 
may choose $g$ and $k$ so that $\sigma$ and 
$\tau$ commute on $C$. By \cite[Proposition 8.8]{ATV2}, $\tau$ induces
an automorphism $\tau'$ of $A$. Twisting $A$ by $\tau'$
produces an algebra $B$ whose point scheme is $C$ with automorphism
$\tau \circ \sigma \in$ Aut$(C)$. Since each line of $C$
is invariant under $\tau \circ \sigma$, it follows, by 
Corollary~\ref{cor1}, that $B$ is a twist, by an automorphism, of an Ore 
extension $D$ of the polynomial ring on two variables.
Hence, $A$ is a twist by an automorphism of a twist by an automorphism 
of $D$. Since twisting by an automorphism need not be transitive, we
can at most conclude that $A$ is a twist of $D$ by a twisting system.

Instead, suppose $b^2 + b + 1 = 0$. If char$(\kk)\neq 3$, then the Koszul 
dual of $A$ has Hilbert series $H(t) = 1 + 3t + 3t^2$, giving that $A$ is not
a quantum $\PP^2$, which is a contradiction. However, if char$(\kk)= 3$, then
$b = 1$, and $A$ is generated by $x,\ y,\ z$ with defining relations:
\[
\begin{array}{l}
x y = x^2 + y^2,\\[2mm]
z y = - x z + c x^2 + e y^2,\\[2mm]
z x = (y -x) z -a y x + c x^2.
\end{array}
\tag{$*$}
\]
In this case, $A = B[z; \phi,\ \delta]$ where $B = \kk\langle x, \
y \rangle/\langle x^2 + y^2 - x y\rangle$ and $\phi \in$ Aut$(B)$
is given by $\phi(x) = y - x$, $\phi(y) = - x$,  and $\delta$
is the left $\phi$-derivation of $B$ given by 
$\delta(x) = cx^2 -a y x$, 
$\delta(y) = cx^2 +e y^2$. 
Mapping $x \mapsto r_2$ and $y \mapsto r_1 - r_2$ yields that $B \cong 
\kk\langle r_1, \ r_2 \rangle/\langle r_1r_2 - r_2r_1 - r_1^2 \rangle$
(since char$(\kk) = 3$), which
is a regular algebra of global dimension two, so, by \cite[Theorem~4.2]{E},
any algebra with the relations \thetag{$*$} is regular.
\end{pf}

Summarizing our work in this section yields the next result.

\begin{thm}\label{main1}
If the point scheme of a quadratic quantum $\PP^2$ is reducible or 
non-reduced, then the algebra is either \/ {\rm (a)}
an Ore extension of a regular algebra of global dimension two, or\/ {\rm
(b)}
a twist, by an automorphism, of a skew polynomial ring, or\/ {\rm (c)} a 
twist, by a 
twisting system, of an Ore extension of the polynomial ring on two variables;
if\/ $\char(\kk) \neq 3$, then the algebra is described by\/ {\rm (b)} or\/
{\rm (c)}.
\end{thm}
\begin{pf}
Combine Corollary~\ref{cor1}, Lemma~\ref{lemma1} and
Proposition~\ref{prop2}.
\end{pf}

If char$(\kk) \neq 2$, then skew polynomial rings are graded skew Clifford 
algebras, so, in this setting, Theorem~\ref{main1} may be rephrased as 
Theorem~\ref{main1'} below.  We first recall the definition
of a graded skew Clifford algebra and of some terms used in its definition.
For the definition, we assume char$(\kk) \neq 2$.

\begin{defn}\cite{CV}\label{GSCA}
For $\{i, j \}\subset \{ 1, \ldots , n\}$, let $\mu_{ij} \in \kk^\times$ 
satisfy $\mu_{ij} \mu_{ji} = 1$ for all $i \neq j$, and 
write $\mu = (\mu_{ij}) \in M(n, \kk)$.  
A matrix $M \in M(n, \kk)$ is called $\mu$-{\em symmetric} if
$M_{ij} = \mu_{ij}M_{ji}$ for all $i, j = 1, \ldots , n$. Henceforth, suppose
$\mu_{ii} = 1$ for all $i$, and fix $\mu$-symmetric matrices $M_1, \ldots , 
M_n \in M(n, \kk)$. 
A {\em graded skew Clifford algebra} 
associated to $\mu$ and $M_1$, $\ldots ,$ $M_n$ is a graded $\kk$-algebra
on degree-one generators $x_1, \ldots , x_n$ and on degree-two generators
$y_1, \ldots , y_n$ with defining relations given by:
\begin{enumerate}
\item[(a)] $\ds x_i x_j + \mu_{ij} x_j x_i = \sum_{k=1}^n (M_k)_{ij} y_k$
           for all $i, j = 1, \ldots , n$, and
\item[(b)] the existence of a normalizing sequence $\{ r_1, \ldots , r_n\}$
           of homogeneous elements of degree two that span 
	   $\kk y_1 + \cdots + \kk y_n$. 
\end{enumerate}
\end{defn}
\noindent 
One should note that if $\mu_{ij} = 1$ for all $i, j$, and if the $y_k$ are 
all central in the algebra, then the algebra is a graded Clifford algebra. 
Moreover, polynomial rings, and skew polynomial rings, on finitely-many 
generators are graded skew Clifford algebras.
Although graded skew Clifford algebras need not be quadratic nor regular
in general, a simple geometric criterion was established in 
\cite[Theorem~4.2]{CV} for determining when such an algebra is quadratic and 
regular.  We refer the reader to \cite{CV,StV2} for results on graded
Clifford algebras and graded skew Clifford algebras.

\begin{ex}
Suppose that char$(\kk) \neq 2$ and that the $\mu_{ij}$ are given as in
Definition~\ref{GSCA}. Fix $\alpha_1,\ \alpha_2,\ \alpha_3 \in \kk$ and 
define
\[
M_1 = \begin{bmatrix}2&0&0\\0&0&\alpha_1\\0&\mu_{32}\alpha_1&0\end{bmatrix},
\quad
M_2 = \begin{bmatrix}0&0&\alpha_2\\0&2&0\\ \mu_{31}\alpha_2&0&0\end{bmatrix},
\quad
M_3 = \begin{bmatrix}0&\alpha_3&0\\ \mu_{21}\alpha_3&0&0\\0&0&2\end{bmatrix}.
\]
These $\mu$-symmetric matrices, with values for $\alpha_1,\ \alpha_2,\ 
\alpha_3$ and the 
$\mu_{ij}$ given below, yield a regular graded skew Clifford algebra of 
global dimension three on generators $x_1,\ x_2,\ x_3$ with defining
relations 
\[
x_i x_j + \mu_{ij} x_j x_i = \alpha_k x_k^2,
\]
for all distinct $i, j, k$,
with point scheme isomorphic to a subscheme $\mathcal P$ (given below)
of $\PP^2$.
\begin{enumerate}
\item[(i)] $\alpha_i = 0$ for all $i$, $\mu_{13} + \mu_{12}\mu_{23} = 0$, 
           $\mathcal P = \PP^2$;
\item[(ii)] $\alpha_i = 0$ for all $i$, $\mu_{13} + \mu_{12}\mu_{23} \neq 0$, 
           $\mathcal P = \V(x_1 x_2 x_3)$, which is a ``triangle'';
\item[(iii)] $\alpha_1 = 0 = \alpha_2\neq\alpha_3$, $\mu_{13}\mu_{23} = 1$, 
             $\mu_{13} + \mu_{12}\mu_{23} = 0$, $\mathcal P = \V(x_3^3)$, 
	     which is a triple line;
\item[(iv)] $\alpha_1 = 0 = \alpha_2\neq\alpha_3$, $\mu_{13}\mu_{23} = 1$, 
             $\mu_{13} + \mu_{12}\mu_{23} \neq 0$, 
	     $\mathcal P = \V(((\mu_{13} + \mu_{12}\mu_{23}) x_1 x_2 +
	     \alpha_3 x_3^2)x_3)$, which is the union of a nondegenerate
	     conic and a line;
\item[(v)] $\alpha_1 = 0 \neq \alpha_2\alpha_3$, $\mu_{12} = \mu_{23}$, 
             $\mu_{12}^3 = 1$, $\mu_{13} = \mu_{12}^2$, 
	     $\mathcal P = \V(\alpha_2 x_2^3 + \alpha_3 \mu_{13}x_3^3 +
	     2 \mu_{12} x_1 x_2 x_3)$, which is a nodal cubic curve;
\item[(vi)] $\alpha_1\alpha_2\alpha_3\neq 0$, $\mu_{12} = \mu_{23}$, 
             $\mu_{12}^3 = 1$, $\mu_{13} = \mu_{12}^2$,
	     $\alpha_1\alpha_2\alpha_3 + \mu_{12}^2 \neq 0$,
	     $\mathcal P = \V(\alpha_1 x_1^3 + \alpha_2 \mu_{12}x_2^3+
	     \alpha_3 x_3^3 + (2 \mu_{12}^2-\alpha_1\alpha_2\alpha_3)
	     x_1 x_2 x_3)$, which is an elliptic curve if and only if 
	     $\alpha_1\alpha_2\alpha_3 \neq 8\mu_{12}^2$.
\end{enumerate}
The method to find the above values uses \cite[Theorem~4.2]{CV}.
For other values of the $\alpha_i$ that yield a regular algebra, see 
\cite{Nafari}. 
This example highlights the wide variety of point schemes that can be 
obtained directly from regular graded skew Clifford algebras of global 
dimension three.
\end{ex}

\begin{thm}\label{main1'}
Suppose\/ $\char(\kk) \neq 2$. 
If the point scheme of a quadratic quantum $\PP^2$ is reducible or 
non-reduced, then either the algebra is a twist, by an automorphism, of a 
graded skew Clifford algebra, or the algebra is a twist, by a twisting 
system, of an Ore extension of a regular graded skew Clifford algebra of
global dimension two.
\end{thm}
\begin{pf}
This is a restatement of Theorem~\ref{main1} in terms of graded skew
Clifford algebras.
\end{pf}

The subsequent sections of the article focus on the case where the cubic
divisor $C$ is reduced and irreducible.  
It is straightforward to prove that if $C$ contains two or more singular 
points, then $C$ contains a line. Thus, in the remaining sections, $C$
has at most one singular point. Indeed, by Lemma~\ref{curve}, $C$ will be 
either a nodal cubic curve, a cuspidal cubic curve, or an elliptic curve.

\bigskip
\bigskip


\section{Quantum $\PP^2$s with Point Scheme a Nodal Cubic Curve}\label{Sec2}

In this section, we classify those quadratic quantum~$\PP^2$s whose point 
scheme is a nodal cubic curve in $\PP^2$, and prove that, up to isomorphism, 
there is at
most a one-parameter family of such algebras (Theorem~\ref{nodal}).
Moreover, if char$(\kk)\neq 2$, we show, in Corollary~\ref{cor2}, that such 
algebras are Ore extensions of regular graded skew Clifford algebras of 
global dimension two, and, under certain conditions, are even graded skew 
Clifford algebras themselves.

Throughout this section, we use $x, y$ and $z$ for homogeneous degree-one 
linearly independent (commutative) coordinates on $\PP^2$. 
Our next result shows that a nodal cubic curve and a cuspidal cubic curve 
are the only irreducible cubic divisors in $\PP^2$ with a unique singular 
point; for lack of a suitable reference, we include its simple proof.

\begin{lemma} \label{curve}
Let $C$ denote an irreducible cubic divisor in $\PP^2$ with a unique singular 
point.  Up to isomorphism, $C = \V(f)$, where either\/
{\rm (a)} $f = x^3 + y^3 + xyz$, or\/
{\rm (b)} $f = y^3 + x^2 z$, or\/ 
{\rm (c)} $f = y^3 + x^2 z + x y^2$; 
if\/ $\char(\kk) \neq 3$, then $f$ is given by\/ {\rm (a)} or\/ {\rm (b)}.
\end{lemma}
\begin{pf}
By rechoosing $x$, $y$ and $z$ if needed, we may assume that $\V(x, y)$ 
is the unique singular point on $C$ and that $C = \V(f)$, where 
$f = s_1 + s_2 xz$, where $s_1 = \alpha_1 x^3 + \alpha_2 x^2 y + 
\alpha_3 x y^2 + y^3$, $s_2 \in \{ x, y\}$ and $\alpha_i \in \kk$ for all 
$i$. Moreover, if $s_2 = y$, then $\alpha_1 \neq 0$, as $C$ is irreducible.

If $s_2 = y$, then $f \mapsto x^3 + y^3 + xyz$ by mapping $x \mapsto 
\beta^{-1}x$, $y \mapsto y$ and 
$z \mapsto \beta z - \alpha_2 \beta^{-1} x - \alpha_3 y$, where
$\beta\in\kk$ is any solution of the equation $\beta^3 = \alpha_1$.

On the other hand, suppose $s_2 = x$. If char$(\kk) \neq 3$, then 
$f \mapsto y^3 + x^2 z$ by mapping $x \mapsto x$, 
$y \mapsto y - (\frac{\alpha_3}{3})x$  and 
$z \mapsto z + (\frac{\alpha_2 \alpha_3}{3} - \alpha_1 - 
\frac{2 \alpha_3^3}{27})x + (\frac{\alpha_3^2}{3} - \alpha_2)y$.
If char$(\kk) = 3$ and $\alpha_3 = 0$, then 
$f \mapsto y^3 + x^2 z$ by mapping $x \mapsto x$, $y \mapsto y$ and 
$z \mapsto z - \alpha_1 x - \alpha_2 y$.
If char$(\kk) = 3$ and $\alpha_3 \neq 0$, then $f \mapsto 
y^3 + x^2 z + x y^2$ by mapping $x \mapsto \alpha_3^{-1}x$, $y \mapsto y$ and 
$z \mapsto \alpha_3^2 z - \alpha_1\alpha_3^{-1}x - \alpha_2 y$.
\end{pf}

In Lemma~\ref{curve}, if $C$ is given by (a), we refer 
to $C$ as a nodal cubic curve; otherwise, we refer to $C$ as a cuspidal
cubic curve.  For the rest of this section, $C$ denotes a nodal cubic curve.

In order to classify those quadratic quantum $\PP^2$ whose point scheme
is given by a nodal cubic curve $C$ in $\PP^2$, we will classify all
such algebras whose defining relations vanish on the graph of an 
automorphism of $C$ (see Remark~\ref{rmk1}(ii)).

\begin{thm}\label{nodal}
Let $A$ denote a quadratic quantum $\PP^2$.  If the point scheme of $A$ is 
a nodal cubic curve in $\PP^2$, then $A$ is isomorphic to a quadratic 
algebra on three generators $x_1$, $x_2$, $x_3$ with defining relations:
\[ 
\lambda x_1 x_2 = x_2 x_1, \quad 
\lambda x_2 x_3 = x_3 x_2 - x_1^2, \quad 
\lambda x_3 x_1 = x_1 x_3 - x_2^2, 
\]
where $\lambda \in \kk$ and $\lambda (\lambda^3 -1) \neq 0$. Moreover,
for any such $\lambda$, any quadratic algebra with these defining relations 
is a quadratic quantum $\PP^2$ with point scheme given by a 
nodal cubic curve in $\PP^2$.  
\end{thm}
\begin{pf}
By Lemma~\ref{curve}, we may write $C=\V(x^3 + y^3 + xyz)$ for the nodal
cubic curve.  It follows that 
\[C = \{ (a^2, a, -a^3 - 1) : a \in \kk\},\] 
and the unique singular point of $C$ is $p = (0, 0, 1)$. Thus if $\sigma
\in\text{Aut}(C)$, then 
\[\sigma(a^2, a, -a^3 - 1)  = (f(a)^2,\ f(a),\ -f(a)^3 - 1),\] 
for all $a \in \kk$, where $f$ is a rational function of one variable.  
Since $\sigma$ is invertible, $f$
is also, so $f(a) = (\lambda_1 a + \lambda_2)/%
(\lambda_3 a + \lambda_4)$, where $\lambda_i \in \kk$ for all $i$.  However, 
$\sigma(p) = p$ implies that $f(0) = 0$, and the domain of $f$ is $\kk$, so 
$\lambda_2 = 0 = \lambda_3$. Thus, there exists $\lambda \in \kk^\times$ such
that $f(a) = \lambda a$ for all $a \in \kk$. It follows that 
\[\sigma(a^2, a, -a^3 - 1)  = (\lambda^2 a^2,\ \lambda  a,\   
-\lambda^3 a^3 - 1),\] 
for all $a \in \kk$.  Since $\sigma$ may be extended to $\PP^2$ if 
$\lambda^3 = 1$, by Remark~\ref{rmk1}(iii), we may assume $\lambda^3 \neq 1$. 

Let $\{x_1 , x_2 , x_3\}$ be a basis for $A_1$, let $\{z_1 , z_2 , z_3\}$ 
be the dual basis for $A\dual_1$, and let $W$ and $W^\perp$ be as in 
Remark~\ref{rmk1}(i).  We may produce some elements of $W^\perp$ from the 
graph of $\sigma$ as follows.
Firstly, suppose that char$(\kk) \neq 3$ and fix $\omega\in \kk$ such that 
$\omega^2 -\omega + 1 = 0$. This yields the following six elements 
in $W^\perp$ corresponding to the given distinct values of $a$:\\[-4mm]
\[
\begin{array}{rl}
a = 0:&z_3^2 \\[2mm]
a = -1:&(z_1 - z_2)(\lambda^2 z_1 - \lambda z_2 + (\lambda^3-1) z_3)\\[2mm]
a = \omega:&\omega (\omega z_1 + z_2)(\lambda^2 \omega^2 z_1 + \lambda
\omega z_2 + (\lambda^3 -1) z_3)\\[2mm] 
a = -\lambda^{-1}:&%
(\lambda^{-2} z_1 - \lambda^{-1} z_2-(1 - \lambda^{-3}) z_3)(z_1 - z_2)\\[2mm]
a = \omega\lambda^{-1}:&%
(\omega^2\lambda^{-2} z_1 + \omega\lambda^{-1} z_2-(1 - \lambda^{-3}) z_3)%
(\omega z_1 + z_2)\omega\\[2mm]
a = -\omega^2:&-\omega (z_1 + \omega z_2)(-\lambda^2 \omega z_1 - \lambda
\omega^2 z_2 + (\lambda^3 -1) z_3).
\end{array}
\]
\quad\\[-1mm]
Taking linear combinations of these six elements yields the following basis 
for $W^\perp$:
\[
\begin{array}{ll}
z_3^2, & z_1 z_2 + \lambda z_2 z_1,\\[2mm]
\lambda z_2^2 + (\lambda^3 - 1) z_1 z_3, &z_2 z_3 + \lambda z_3 z_2,\\[2mm]
\lambda z_1^2 + (\lambda^3 - 1) z_3 z_2, \qquad & z_3 z_1 + \lambda z_1 z_3.
\end{array}
\]
It follows that $W$ is the span of the elements:
\[ 
\lambda x_1 x_2 - x_2 x_1, \quad 
\lambda ( \lambda x_2 x_3 - x_3 x_2) + (\lambda^3-1) x_1^2, \quad 
\lambda ( \lambda x_3 x_1 - x_1 x_3) + (\lambda^3-1) x_2^2,  
\]
if char$(\kk) \neq 3$.  Moreover,  these three linearly independent
elements vanish on the graph 
of $\sigma$ even if char$(\kk) = 3$,  so $W$ is the span of these three 
elements even in this case.  Furthermore, since $\lambda(\lambda^3 -1) 
\neq 0$, we may map $x_1 \mapsto x_1$, $x_2 \mapsto x_2$  and 
$x_3 \mapsto \lambda^{-1}(\lambda^3 -1)x_3$, so $A$ is isomorphic to the
algebra in the statement of the theorem.

If $\lambda \in \kk$ where $\lambda (\lambda^3 -1) \neq 0$, then
an algebra with the given relations is an Ore extension $B[x_3; \phi,\
\delta]$ of the algebra~$B = 
\kk\langle x_1 , x_2\rangle/\langle \lambda x_1 x_2 - x_2 x_1\rangle$
using $\phi \in \text{Aut}(B)$ and $\delta$ a left $\phi$-derivation of $B$ 
where 
\[ 
\phi(x_1) = \lambda^{-1}x_1, \quad  
\phi(x_2) = \lambda x_2, \quad  
\delta(x_1) = -\lambda^{-1}x_2^2, \quad  
\delta(x_2) = x_1^2.  
\]
Since $B$ is a regular algebra of global dimension two, it follows, by 
\cite[Theorem~4.2]{E}, that such an Ore extension of $B$ is a regular algebra
of global dimension three.

The point scheme of the algebra with the defining relations in the
theorem is given by 
$\V( \lambda x^3 + \lambda y^3 +(\lambda^3 -1) x y z)$, which is 
indeed a nodal cubic curve.  
\end{pf}

\begin{cor}\label{cor2}
Suppose\/ $\char(\kk) \neq 2$.
If $\lambda^3 = -1$, then the algebra in Theorem~\ref{nodal} is a
graded skew Clifford algebra; if $\lambda^3 \notin \{0, 1\}$, then the
algebra is an Ore extension of a regular graded skew Clifford algebra.
\end{cor}
\begin{pf}
Let $S$ denote the quadratic algebra on generators $z_1$, $z_2$ and
$z_3$ with defining relations 
\[
z_1 z_2 + \lambda z_2 z_1=0,\quad
z_2 z_3 + \lambda z_3 z_2=0,\quad
z_3 z_1 + \lambda z_1 z_3=0.
\]
If $\lambda^3 = -1$, then the set 
$X= \{z_3^2, z_2^2 + z_1 z_3, z_1^2 + z_3 z_2\}$ 
is a normalizing sequence in $S$. In the free algebra on $z_1, z_2, z_3$, 
let $Y$ denote the span of the defining relations of $S$, and let $\hat X$ 
denote the span of any preimage of $X$.  The zero locus in 
$\PP^2\times \PP^2$ of $\hat X + Y$ is the empty set. By
\cite[Theorem~4.2]{CV}, since char$(\kk) \neq 2$, it follows that the
Koszul dual of $S/\langle X \rangle$ is a regular graded skew Clifford 
algebra; by construction, this algebra is isomorphic to the algebra in 
Theorem~\ref{nodal}.

If $\lambda(\lambda^3 - 1) \neq 0$, then the proof of Theorem~\ref{nodal} 
shows that the algebra therein is an Ore extension of a regular algebra~$B$ 
of global dimension two, and $B$ is a graded skew Clifford algebra by 
\cite[Corollary 4.3]{CV}.
\end{pf}

\bigskip
\bigskip


\section{Quantum $\PP^2$s with Point Scheme a Cuspidal Cubic Curve}%
\label{Sec3}

In this section, we prove that, up to isomorphism, there is a unique 
quadratic quantum~$\PP^2$ whose point scheme is a cuspidal cubic curve 
in $\PP^2$ (Theorem~\ref{cusp}).  Moreover, this algebra exists if and only 
if char$(\kk) \neq 3$. We also prove in Corollary~\ref{cor3} that, if 
char$(\kk) \neq 2$, then such 
an algebra is an Ore extension of a regular graded skew Clifford algebra
of global dimension two.

As in Section~\ref{Sec2}, we continue to use $x, y$ and $z$ for homogeneous 
degree-one linearly independent (commutative) coordinates on $\PP^2$.  By 
Lemma~\ref{curve}, we may assume that the cuspidal cubic curve is given
by $C = \V(y^3 + x^2 z)$ or 
$C = \V(y^3 + x^2 z + x y^2)$, with the second case occurring only if 
char$(\kk) = 3$.

\begin{thm}\label{cusp}
Let $A$ denote a quadratic quantum $\PP^2$.  If\/ $\char(\kk) = 3$, then
the point scheme of $A$ is not a cuspidal cubic curve in $\PP^2$.
If\/ $\char(\kk) \neq 3$ and if the point scheme of $A$ is a cuspidal cubic 
curve in $\PP^2$, then $A$ is isomorphic to a quadratic algebra on three 
generators $x_1$, $x_2$, $x_3$ with defining relations:
\[ 
x_1 x_2 = x_2 x_1 + x_1^2, \quad 
x_3 x_1 = x_1 x_3 + x_1^2 + 3 x_2^2, \quad 
x_3 x_2 = x_2 x_3 - 3 x_2^2 - 2 x_1 x_3 - 2 x_1 x_2. 
\]
Moreover, any quadratic algebra with these defining relations is a quadratic 
quantum $\PP^2$; it has point scheme given by a cuspidal cubic curve in 
$\PP^2$ if and only if\/ $\char(\kk) \neq 3$. 
\end{thm}
\begin{pf}
Suppose first that the cuspidal cubic curve is $C=\V(y^3 + x^2 z)$.  It 
follows that 
\[C = \{(0, 0, 1)\} \cup \{ (1, b, -b^3) : b \in \kk\},\] 
and that the unique singular point of $C$ is $p = (0, 0, 1)$. 
Thus, if $\sigma \in\text{Aut}(C)$, then $\sigma (p) = p$ and 
\[\sigma(1, b, -b^3)  = (1,\ f(b),\ -f(b)^3),\] 
for all $b \in \kk$, where $f$ is a rational function of one variable.  
Since $\sigma$ is invertible, $f$ is also, so 
$f(b) = (\lambda_1 b + \lambda_2)/%
(\lambda_3 b + \lambda_4)$, where $\lambda_i \in \kk$ for all $i$.  
However, the domain of $f$ is $\kk$, so $\lambda_3 = 0$ and $\lambda_1,
\lambda_4 \in \kk^\times$.   Writing the points of $C$ in the form 
$(a^3, a^2 b, -b^3)$ for all $(a, b) \in \PP^1$ and rechoosing the 
$\lambda_i$, we may write 
\[ \sigma(a^3, a^2 b, -b^3) = (a^3,\ \lambda_1 a^2 (b + \lambda_2 a),\  
-\lambda_1^3 (b + \lambda_2 a)^3)
\]
for all $(a, b) \in \PP^1$, where $\lambda_1 \in \kk^\times$ and
$\lambda_2 \in \kk$.  If $\lambda_2 = 0$ or if char$(\kk) = 3$, then $\sigma$
may be extended to $\PP^2$, so, by Remark~\ref{rmk1}(iii), we may assume  
$\lambda_2 \neq 0$ and char$(\kk) \neq 3$. 

Using the method and notation in the proof of Theorem~\ref{cusp}, we find 
that $W^\perp$ has basis:
\[
\begin{array}{ll}
z_3^2, & z_2^2 - 3 \lambda_1^2 \lambda_2 z_1 z_3 - 3 \lambda_1^2
\lambda_2^2 z_2 z_3,\\[2mm]
z_3 z_2 + \lambda_1^2 z_2 z_3,&
z_3 z_1 + \lambda_1^3 z_1 z_3 + 2 \lambda_1^3 \lambda_2 z_2 z_3,\\[2mm]
z_1^2 + \lambda_1 \lambda_2 z_1 z_2 - \lambda_1^3 \lambda_2^3 z_1 z_3 , 
\qquad & 
z_2 z_1 + \lambda_1 z_1 z_2 + 2 \lambda_1^3 \lambda_2^3 z_2 z_3.
\end{array}
\]
(Alternatively, the reader may simply verify that the dual elements to
these elements vanish on the graph of $\sigma$.)
Mapping $z_1 \mapsto z_1$, 
$z_2 \mapsto z_2/\lambda_2$ and
$z_3 \mapsto z_3/\lambda_2^3$ allows us to take $\lambda_2 = 1$.
It follows that the Hilbert series of the Koszul dual of $A$ is 
$H(t) = (1+t)^3$ if and only if $\lambda_1=1$ (since char$(\kk) \neq 3$).
If $\lambda_1=1$,
then $A$ is the algebra given in the statement of the theorem, where
$\{x_1 , x_2 , x_3\}$ is the dual basis to $\{z_1 , z_2 , z_3\}$.

To prove the relations in the statement determine a regular algebra, we 
write the algebra as an Ore extension of the regular 
algebra~$B = \kk\langle x_1 , x_2\rangle/\langle x_1 x_2 - x_2 x_1 - 
x_1^2\rangle$ using $\phi \in \text{Aut}(B)$ and $\delta$ a left 
$\phi$-derivation of $B$ where 
\[ 
\phi(x_1) = x_1, \quad  
\phi(x_2) = x_2 - 2 x_1, \quad  
\delta(x_1) = x_1^2 + 3 x_2^2, \quad  
\delta(x_2) = -2 x_1 x_2 - 3 x_2^2.  
\]
It follows, from \cite[Theorem~4.2]{E}, that such an Ore extension of~$B$ 
is a regular algebra of global dimension three.  
The point scheme of such an Ore
extension is given by $\V( 3 (y^3 + x^2 z) )$, which is a
cuspidal cubic curve if and only if char$(\kk) \neq 3$.

By Lemma~\ref{curve}, we henceforth assume $C = \V(y^3 + x^2 z + x y^2)$ 
and char$(\kk) = 3$.  It follows that
\[C = \{(0, 0, 1)\} \cup \{ (1, b, -b^2-b^3) : b \in \kk\},\] 
and that the unique singular point of $C$ is $p = (0, 0, 1)$. 
In this setting, if $\sigma \in\text{Aut}(C)$, then $\sigma (p) = p$ and 
\[\sigma(1, b, -b^2-b^3)  = (1,\ f(b),\ -f(b)^2-f(b)^3),\] 
for all $b \in \kk$, where $f$ is a rational function of one variable.  
As before, the invertibility of $\sigma$ implies that $f$ is invertible,
and that $f(b) = (\lambda_1 b + \lambda_2)/%
(\lambda_3 b + \lambda_4)$, where $\lambda_i \in \kk$ for all $i$.  
Since the domain of $f$ is $\kk$, $\lambda_3 = 0$ and $\lambda_1, \lambda_4 
\in \kk^\times$. Thus, writing the points of $C$ in the form 
$(a^3, a^2 b, -ab^2 -b^3)$ for all $(a, b) \in \PP^1$ and rechoosing the 
$\lambda_i$, we may write 
\[ \sigma(a^3, a^2 b, -ab^2 -b^3) = (a^3,\ a^2 (\lambda_1 b + \lambda_2 a),\  
- a(\lambda_1 b + \lambda_2 a)^2 - (\lambda_1 b + \lambda_2 a)^3)
\]
for all $(a, b) \in \PP^1$, where $\lambda_1 \in \kk^\times$ and
$\lambda_2 \in \kk$.  By Remark~\ref{rmk1}(iii), we further assume that 
$\lambda_1 \neq 1$, since $\lambda_1 = 1$ if and only if $\sigma$ may be
extended to $\PP^2$ (since char$(\kk) = 3$).

Using the notation in the proof of Theorem~\ref{cusp},  and using 
$x_1, x_2 , x_3$ as generators for $A$, and by seeking homogeneous degree-two
elements that vanish on the graph of $\sigma$, we find that a basis for $W$ 
is:\\[-5mm]
\[
\begin{array}{l}
x_1 x_2 - \lambda_1 x_2 x_1 - \lambda_2 x_1^2,\\[2mm]
x_1 x_3 - \lambda_1^3 x_3 x_1 + \lambda_1 (1 - \lambda_1)x_2^2 +
\lambda_1 \lambda_2 (1 + \lambda_1) x_2 x_1 + \lambda_2^2 (1 + \lambda_2) 
x_1^2,\\[2mm]
x_2 x_3 - \lambda_1^2 x_3 x_2 + \lambda_1^2 (\lambda_1 + \lambda_2 -1)
x_3 x_1 + (\lambda_1^2 -\lambda_1 + 2\lambda_2) x_2^2 + \lambda_2(\lambda_2^2
- \lambda_2 - \lambda_1^2 + \lambda_1 ) x_2 x_1.
\end{array}
\]
\\[-2mm]
Since $\lambda_1(\lambda_1-1) \neq 0$, it follows that the Hilbert series 
of $A$ is $H(t) = 1 + 3t + 6t^2 + 9t^3 + \cdots$, which contradicts $A$ 
being a quadratic quantum $\PP^2$.
\end{pf}

\begin{cor}\label{cor3}
If\/ $\char(\kk) \neq 2$, then the algebra in Theorem~\ref{cusp} is an Ore 
extension of a regular graded skew Clifford algebra.
\end{cor}
\begin{pf}
The proof of Theorem~\ref{cusp} shows that the algebra is an Ore extension 
of a regular algebra~$B$ of global dimension two, and $B$ is a graded skew 
Clifford algebra by \cite[Corollary~4.3]{CV}.
\end{pf}

\bigskip
\bigskip


\section{Quantum $\PP^2$s with Point Scheme an Elliptic Curve}\label{Sec4}

It remains to consider quadratic quantum $\PP^2$s with point scheme an 
elliptic curve. In \cite{AS}, such algebras are classified into four types, 
A, B, E and H, where some members of each type might not have an elliptic 
curve as their point scheme, but a generic member does.
We show, in Lemmas~\ref{H}, \ref{B} and \ref{A}, 
that all regular algebras of types B and H that have point scheme an 
elliptic curve, and some regular algebras of type A that have point scheme 
an elliptic 
curve are given by graded skew Clifford algebras, or twists thereof.
Up to isomorphism and anti-isomorphism, type~E consists of at most one 
algebra and it appears not to be directly related to a graded skew Clifford 
algebra (but this issue is still open), so this type is only discussed in
Remark~\ref{weak} regarding a weak relationship to a graded skew Clifford 
algebra.

\subsection{Type H}\hfill

By \cite[Page 207]{AS}, there are at most two regular algebras of type H 
(up to isomorphism) and they are given by 
$\kk$-algebras on generators $x, y, z$ with defining relations:
\[
y^2 = x^2, \quad zy = - i yz, \quad  yx - xy = i z^2,
\]
where $i$ is a primitive fourth root of unity.  In the following result, such
an algebra is denoted~$H$; its point scheme is an elliptic curve unless 
char$(\kk) = 2$.

\begin{lemma}\label{H}
If\/ $\char(\kk) \neq 2$, then the algebra $H$ is a regular graded skew 
Clifford algebra and a twist of a graded Clifford algebra by an automorphism.
\end{lemma}
\begin{pf}
Suppose char($\kk) \neq 2$, and let $H\dual_1$ have basis $\{X, Y, Z\}$ 
dual to $\{x, y, z\}$.   Let $S$ denote the $\kk$-algebra on $X, Y, Z$
with defining relations:
\[
YX = - XY, \qquad  YZ = i ZY, \qquad ZX = \nu XZ,
\]
where $\nu \in \kk^\times$. 
For all $\nu \in \kk^\times$, the Koszul dual $H^{\sstar}$ to $H$ is the 
quotient of $S$ by the ideal spanned 
by the normalizing sequence $\{ XZ,\ i XY - Z^2,\ X^2 + Y^2\}$.
The defining relations of $H^{\sstar}$ have empty zero locus in $\PP^2 \times 
\PP^2$. By \cite[Theorem~4.2]{CV},
$H$ is a regular graded skew Clifford algebra. If, 
further, we choose $\nu = i$, then $S$ is a twist of a polynomial ring
by an automorphism,
so, by \cite[Proposition 4.5]{CV}, $H$ is a twist of a graded Clifford
algebra by an automorphism.
\end{pf}

\subsection{Type B}\quad\\
\indent 
By \cite[Page 207]{AS}, the regular algebras of type B that have point scheme
an elliptic curve are given by $\kk$-algebras on generators $x, y, z$ with 
defining relations:
\[
xy + yx = z^2 - y^2, \quad 
xy + yx = az^2 - x^2, \quad 
zx - xz = a(yz -zy),
\]
where $a \in \kk$, $a(a-1) \neq 0$. (A sign error in the first relation on 
Page~207 of \cite{AS} has been corrected above.) In the 
following result, such an algebra is denoted~$B$; its point scheme is an 
elliptic curve for generic values of $a$ unless char$(\kk) \in \{ 2, \ 3\}$.

\begin{lemma}\label{B}
If\/ $\char(\kk) \neq 2$, then the algebra $B$ is regular if and only if 
$a^2 -a + 1 \neq 0$; in this case, $B$ is a graded skew Clifford algebra 
and a twist of a graded Clifford algebra by an automorphism.
\end{lemma}
\begin{pf}
Suppose char($\kk) \neq 2$, and let $B\dual_1$ have basis $\{X, Y, Z\}$ 
dual to $\{x, y, z\}$.   Let $S$ denote the $\kk$-algebra on $X, Y, Z$
with defining relations:
\[
YX = XY, \qquad  YZ = - ZY, \qquad ZX = - XZ.
\]
The Koszul dual $B^{\sstar}$ to $B$ is the quotient of $S$ by the ideal 
spanned by the normalizing sequence $\{ Z(a X - Y),\ X^2 + Y^2 -XY,\ 
a X^2 + Y^2 + Z^2\}$.
The defining relations of $B^{\sstar}$ have empty zero locus in $\PP^2 \times 
\PP^2$ if and only if $a^2 -a + 1 \neq 0$. In fact, if $a^2 -a + 1 = 0$, 
then $\dim_\kk(B^{\sstar})$ is infinite, and so $B$ is not regular.
On the other hand, suppose $a^2 -a + 1 \neq 0$. It follows from 
\cite[Theorem~4.2]{CV} that $B$ is a regular graded skew Clifford algebra, 
and, by \cite[Proposition 4.5]{CV}, that $B$ is a twist of a graded Clifford 
algebra by an automorphism, since $S$ is a twist of a polynomial ring by an 
automorphism.
\end{pf}

\subsection{Type A}\hfill

By \cite[Page 207]{AS}, the regular algebras of type A that have point scheme
an elliptic curve are given by $\kk$-algebras on generators $x, y, z$ with 
defining relations:
\[
a xy + b yx + c z^2 = 0, \quad 
a yz + b zy + c x^2 = 0, \quad 
a zx + b xz + c y^2 = 0, \quad 
\tag{$*$}
\]
where $a, b, c \in \kk$.  We denote such an algebra by~$A'$; by \cite{ATV1}, 
its point scheme is an elliptic curve if and only if $abc \neq 0$, 
$(3abc)^3 \neq (a^3 + b^3 + c^3)^3$ and $\char(\kk) \neq 3$.
Thus, we assume $abc \neq 0$ and $(3abc)^3 \neq (a^3 + b^3 + c^3)^3$ and,
with these assumptions, $A'$ is regular unless $a^3 = b^3 = c^3$~(\cite{ATV1}).

\begin{lemma}\label{A}
Suppose\/ $\char(\kk) \neq 2$. 
If $a^3 = b^3 \neq c^3$, then $A'$ is a regular 
graded skew Clifford algebra and a twist of a graded Clifford algebra by an 
automorphism. If $b^3 = c^3 \neq a^3$ or if $a^3 = c^3 \neq b^3$, then
$A'$ is a twist of a regular graded skew Clifford algebra by an automorphism.
\end{lemma}
\begin{pf}
Define $\tau \in$ Aut$(A')$ by $\tau(x) = y$, $\tau(y) = z$ and $\tau(z) = x$.
Twisting $A'$ by $\tau$ (respectively, by $\tau^2$) yields an algebra on 
$x, y, z$ with the same defining relations as in \thetag{$*$} except that 
$a, b$ and $c$ have been cyclically permuted one place (respectively, two
places) to the left. Thus, the second part of the result follows from
the first part, so it remains to prove the first part. 

Let $(A'_1)\dual\ $ have basis $\{X, Y, Z\}$ dual to $\{x, y, z\}$ and let 
$S$ denote the $\kk$-algebra on $X, Y, Z$ with defining relations:
\[
aYX = bXY, \qquad  aZY = bYZ, \qquad aXZ = bZX.
\]
Suppose $a^3 = b^3$. In this case, 
$\{ cXY - a Z^2, \ cYZ - a X^2, \ cZX - a Y^2 \}$ is a 
normalizing sequence in $S$, and the Koszul dual $(A')^{\sstar}$ to $A'$ is 
the quotient of $S$ by the ideal spanned by this normalizing sequence.
The defining relations of $(A')^{\sstar}$ have empty zero locus in $\PP^2 
\times \PP^2$ if $a^3 \neq c^3$.  By \cite[Theorem~4.2]{CV}, it follows that 
$A'$ is a regular graded skew Clifford algebra if $a^3 \neq c^3$. 
By \cite[Proposition 4.5]{CV}, $B$ is a twist 
of a graded Clifford algebra, since $S$ is a twist of a polynomial ring by 
an automorphism (since $a^3 = b^3$).  
\end{pf}

If $abc \neq 0$ and $(3abc)^3 \neq (a^3 + b^3 + c^3)^3$
and $a^3 \neq b^3\neq c^3 \neq a^3$, then it is still open whether or not 
$A'$ is directly related to a graded skew Clifford algebra.

\begin{rmk}\label{weak}
If $\tilde A$ is an algebra of type $A$ or $E$, then the Koszul dual of
$\tilde A$ is the quotient of a regular graded skew Clifford algebra $S$
(indeed, $S$ is a skew polynomial ring). So, in this sense, such algebras 
are weakly related to graded skew Clifford algebras.
\end{rmk}

\begin{ques}
Can the results of this article be generalized to quadratic regular algebras 
of global dimension four, thereby possibly enabling the classification
of such algebras by using regular graded skew Clifford algebras? 
\end{ques}



\bigskip
\bigskip

\begin{thebibliography}{99}
\raggedbottom
\itemsep8pt
{\small
%
\bibitem{A}
{\sc M.~Artin}, Geometry of Quantum Planes, {\it in} ``Azumaya Algebras,
Actions and Modules,'' Eds.~D.~Haile and J.~Osterburg, {\it Contemporary
Math.} {\bf 124} (1992), 1-15.

\bibitem{AS}
{\sc M.~Artin and W.~Schelter}, Graded Algebras of Global
Dimension 3, {\it Adv.\ Math.} {\bf 66} (1987), 171-216.

\bibitem{ATV1}
{\sc M.~Artin,  J.~Tate and M.~Van den Bergh},  Some Algebras
Associated to Automorphisms of Elliptic Curves, {\it The Grothendieck
Festschrift} {\bf 1}, 33-85,
Eds.\ P.~Cartier et al., Birkh\"auser (1990).

\bibitem{ATV2}
{\sc M.~Artin, J.~Tate and M.~Van den Bergh}, Modules over Regular
Algebras of Dimension~3, {\it Invent.\ Math.} {\bf 106} (1991), 335-388.

\bibitem{CV}
{\sc T.\ Cassidy and M.\ Vancliff}, Generalizations of Graded Clifford
Algebras and of Complete Intersections, {\em J.\ Lond.\ Math.\ Soc.}\ 
{\bf 81} (2010), 91-112.

\bibitem{E}
{\sc E.\ K.\ Ekstr\"om}, The Auslander Condition on Graded and Filtered
Noetherian Rings, in: S\'eminaire Dubreil-Malliavin 1987-1988, Lecture
Notes in Mathematics {\bf 1404} (Springer, Berlin, 1989).

\bibitem{Nafari}
{\sc M.\ Nafari}, {\em Regular Algebras Related to Regular Graded Skew 
Clifford Algebras of Low Global Dimension}, Ph.D.\ Thesis, University of 
Texas at Arlington, August 2011.

\bibitem{StV2}
{\sc D.\ R.\ Stephenson and M.\ Vancliff}, Constructing Clifford Quantum
$\PP^3$s with Finitely Many Points, {\em J.~Algebra} {\bf 312} No.~1 (2007),
86-110.

\bibitem{VV}
{\sc M.~Vancliff and K.~Van Rompay}, Embedding a Quantum Nonsingular
Quadric in a Quantum~$\PP^3$, {\it J.~Algebra} {\bf 195} No.~1 (1997),
93-129.

\bibitem{Z}
{\sc J.~J.~Zhang}, Twisted Graded Algebras and Equivalences of Graded 
Categories, {\it Proc.\ London Math.\ Soc.} {\bf 72} No.\ 2 (1996), 281-311.
%
}
%
\end{thebibliography}
\end{document}